\documentclass[12pt]{article}

\usepackage{amsmath,amssymb,amsthm}

\title{A remark on  Chapple-Euler Theorem in non Euclidean geometry}

\author{Takeo Noda and Shin-ichi Yasutomi}

\newtheorem{thm}{Theorem}[section]
\newtheorem{cor}[thm]{Corollary}

\theoremstyle{definition}

\theoremstyle{remark}
\newtheorem{rem}{Remark}[section]

\numberwithin{equation}{section}

\usepackage{color}

\begin{document}

\maketitle

\begin{abstract}
In non-Euclidean geometry, there are several known
correspondings
to Chapple-Euler Theorem.
This remark shows that those results yield  expressions corredponding
to the well-known formula $d=\sqrt{R(R-2r)}$.
\end{abstract}

\section{Chapple-Euler Theorem  in non Euclidean geometry}
For a triangle, let $R$, $r$ and $d$ be its circumradius, inradius and
the distance between its circumcenter and its incenter respectively.
In 1746,
 Chapple  gave  following Theorem in  Euclidean geometry.
\begin{thm}[\cite{Ch}]\label{t1}
\begin{align*}
d=\sqrt{R(R-2r)}.
\end{align*}
\end{thm}

In 1765,
 Euler (See \cite{B}) also reached
the same result,
which is called Chapple-Euler Theorem.
Theorem \ref{t1} implies an inequality $R\geq 2r$, which is called Euler  inequality
or Chapple-Euler inequality.

We consider spherical geometry by a surface of constant curvature $K=1$ and hyperbolic geometry
with
curvature $K=-1$.
In spherical geometry and hyperbolic geometry
similar formulas were given by Cho and Narajo  as follows.
\begin{thm}[\cite{CN}]\label{t2}
Given a triangle in spherical, or hyperbolic geometry following formula holds respectively.
\[\begin{array}{ll}
\sin^2 d=\sin^2 (R-r)-\sin^2 r \cos^2 R &  (\text{spehrical}),\\
\sinh^2 d=\sinh^2 (R-r)-\sinh^2 r \cosh^2 R &  (\text{hyperbolic}).
\end{array}\]
\end{thm}

Alabdullatif  also showed
the
following formula in hyperbolic geometry.

\begin{thm}[\cite{A}]\label{t3}
\begin{align*}
\tanh r=\dfrac{\tanh (R+d)(\cosh^2 R(\sinh^2r-1)+\cosh^2(r+d))}{\cosh^2(r+d)-\cosh^2 R\cosh^2r}.
\end{align*}
\end{thm}

In this article, we give the following Euler-Chapple type formulas in spherical and
hyperbolic geometry, which are equivalent to
the above theorems.

\begin{thm}\label{t4}
Given a triangle in spherical, or hyperbolic geometry following formula holds for in spherical and
hyperbolic geometry, respectively.
\[
\begin{array}{ll}
\tan d=\begin{cases}
\dfrac{\sqrt{\tan R(\tan R-2\tan r)}}{\sqrt{\tan^2 r+(1+\tan r\tan R)^2}},&\text{if $R<\frac{\pi}{2}$,}
\vspace{1mm}\\
0,&\text{if $R=\frac{\pi}{2}$},
\end{cases}
&(\text{spherical}),
\vspace{2mm}\\
\tanh d=\dfrac{\sqrt{\tanh R(\tanh R-2\tanh r)}}{\sqrt{-\tanh^2 r+(1-\tanh r\tanh R)^2}}
&(\text{hyperbolic}).
\end{array}
\]
\end{thm}

This theorem implies the following Euler-Chapple type inequalities in spherical and hyperbolic geometry:
\begin{cor}
The circumradius $R$ and inradius $r$ of a triangle in spherical or hyperbolic geometry satisfy
the following inequality, respectively:
\[
\begin{array}{ll}
\tan R\geq 2\tan r &  (\text{spherical} ),\\
\tanh R\geq 2\tanh r & (\text{hyperbolic} ).
\end{array}
\]
\end{cor}

\begin{rem}
These inequalities have been shown directly by
Svrtan and Veljan \cite{SV}.
\end{rem}

\bigskip
\noindent
{\bf Proof of Theorem \ref{t4}}.

First, we consider a triangle in spherical geometry.
In the case $R=\frac{\pi}{2}$, both of circumscribed and inscribed circles coincide
to a great circle, so $d=0$ or $d=\pi$ and therefore $\tan d=0$.

Thus we may assume $r<R<\frac{\pi}{2}$ in the following.
By Theorem
\ref{t2}, we have
\begin{align*}
\sin^2 d=\sin^2 (R-r)-\sin^2 r \cos^2 R.
\end{align*}
Hence, we have
\begin{align*}
\dfrac{\tan^2d}{1+\tan^2d}=\dfrac{\tan^2(R-r)}{1+\tan^2(R-r)}-\dfrac{\tan^2r}{1+\tan^2r} \dfrac{1}{1+\tan^2R},
\end{align*}
which implies
\begin{align}\label{tan2d}
\tan^2 d(\tan^2 r+(1+\tan r\tan R)^2)=\tan R(\tan R-2\tan r).
\end{align}
From (\ref{tan2d}) and the facts that
$0\leq d<\frac{\pi}{2}$ and $0< r<\frac{\pi}{2}$,
we have Theorem.

Next, we consider a triangle in hyperbolic geometry.
By Theorem
\ref{t2},
we have
\begin{align*}
\sinh^2 d=\sinh^2 (R-r)-\sinh^2 r \cosh^2 R.
\end{align*}
Similarly, we have
\begin{align}\label{tanh2d}
\tanh^2 d(-\tanh^2 r+(1-\tanh r\tanh R)^2)=\tanh R(\tanh R-2\tanh r).
\end{align}
First, we suppose
\begin{align}\label{mtanh2r}
-\tanh^2 r+(1-\tanh r\tanh R)^2=0.
\end{align}
Then, from (\ref{tanh2d}) we have
\begin{align}\label{tanhR}
\tanh R-2\tanh r=0.
\end{align}
From (\ref{mtanh2r}) and (\ref{tanhR}) we have  $\tanh R=1,-1,2,-2$, which contradicts $0<\tanh R<1$.
Next, we suppose
\begin{align}\label{mtanh2r2}
-\tanh^2 r+(1-\tanh r\tanh R)^2<0.
\end{align}
By making the triangle continuously smaller and smaller, we see that $R$ and $r$ are sufficiently small
so that $-\tanh^2 r+(1-\tanh r\tanh R)^2>0$.
Hence, there exists a triangle such that $-\tanh^2 r+(1-\tanh r\tanh R)^2=0$, which
contradicts the  previous consideration.
Therefore, we have
\begin{align}\label{mtanh2r3}
-\tanh^2 r+(1-\tanh r\tanh R)^2>0.
\end{align}
From (\ref{tanh2d}) and (\ref{mtanh2r3}) we have Theorem.
We note that we can also prove the hyperbolic case from Theorem \ref{t3}.
\hfill $\Box$

\vspace{2cm}

\noindent
Takeo Noda: Faculty of Science, Toho University, JAPAN\\
{\it E-mail address: noda@c.sci.toho-u.ac.jp}

\noindent
Shin-ichi Yasutomi: Faculty of Science, Toho University, JAPAN\\
{\it E-mail address: shinichi.yasutomi@sci.toho-u.ac.jp}
\end{document}